\documentclass[xypic,amssymb,rawfonts,amsthm,amsmath,amsfonts,syntonly,11pt]{amsart}
\usepackage{amssymb}

\newcounter{thecounter}
\numberwithin{thecounter}{section}
\newtheorem{lemma}[thecounter]{Lemma}

\newtheorem{thm}[thecounter]{Theorem}

\theoremstyle{definition}

\newtheorem{rem}[thecounter]{Remark}

\numberwithin{equation}{section}

\newcommand{\Tor}{\operatorname{Tor}}

\newcommand{\pcom}{\hat{{}_p}}

\newcommand{\twocom}{\hat{{}_2}}

\newcommand{\Sp}{\operatorname{Sp}}

\newcommand{\Lie}{\operatorname{Lie}}

\newcommand{\Q}{{\mathbf {Q}}}

\newcommand{\F}{{\mathbf {F}}}
\newcommand{\Z}{{\mathbf {Z}}}

\newcommand{\A}{{\mathbf A}}

\newcommand{\C}{{\mathbf{C}}}

\newcommand{\beq}{\begin{eqnarray*}}
\newcommand{\eeq}{\end{eqnarray*}}
\newcommand{\pil}[1]{\stackrel{#1}{\rightarrow}}

\newcommand{\tuborg}{\left\{\begin{array}{ll}}
\newcommand{\sluttuborg}{\end{array}\right.}

\newfont{\bm}{msbm10}

\newcommand{\Lin}{\operatorname{Lin}}

\newcommand{\dT}{\breve T}

\begin{document}
\title{A finite loop space not rationally equivalent to a compact Lie group} 


\author[K. Andersen]{Kasper K. S. Andersen}
\author[T. Bauer]{Tilman Bauer}
\author[J. Grodal]{Jesper Grodal}
\author[E. K. Pedersen]{Erik Kj\ae r Pedersen}

\thanks{The third named author was partially supported by NSF grant DMS-0104318}
\thanks{The fourth named author thanks SFB 478 M\"unster for its
  hospitality and support}

\subjclass{Primary: 55P35; Secondary: 55P15, 55R35}

\address{Department of Mathematics\newline
\indent University of Copenhagen, Universitetsparken 5\newline
\indent DK-2100 Copenhagen \O, Denmark}
\email{kksa@math.ku.dk}
\address{Sonderforschungsbereich 478 Geometrische Strukturen in der
Mathematik\newline
\indent Westf\"alische Wilhelms-Universit\"at M\"unster, Hittorfstr. 
27\newline
\indent D-48149 M\"unster, Germany}
\email{tbauer@math.uni-muenster.de}
\address{Department of Mathematics\newline
\indent University of Chicago, 5734 S. University Ave.\newline
\indent Chicago, IL 60637, USA}
\email{jg@math.uchicago.edu}
\address{Dept. of Mathematical Sciences\newline
\indent Binghamton University\newline
\indent Binghamton, New York 13902-6000, USA}
\email{erik@math.binghamton.edu}

\begin{abstract}
We construct a connected finite loop space of rank $66$ and dimension
$1254$ whose rational cohomology is not isomorphic as a graded vector
space to the rational cohomology of any compact Lie group, hence
providing a counterexample to a classical conjecture. Aided by machine
calculation we verify that our counterexample is minimal, i.e., that
any finite loop space of rank less than $66$ is in fact rationally
equivalent to a compact Lie group, extending the classical known bound
of $5$.
\end{abstract}

\maketitle

\section{Introduction}
Since the discovery of the Hilton-Roitberg `criminal'
\cite{HR68,HR69,stasheff69} in $1968$ it has been clear that not every
finite loop space is homotopy equivalent to a compact Lie group. The
conjecture emerged however, that this should hold rationally, i.e.,
that any finite loop space should be rationally equivalent to some
compact Lie group (see \cite{AW80}, \cite[~p.~67]{kane88}), and
evidence for this has been accumulated over the years.
In this paper we resolve this conjecture in the negative by
exhibiting a concrete finite loop space of rank $66$ whose rational
cohomology does not agree with that of any compact Lie group.  To do this, we
first use Sullivan's arithmetic square
\cite{sullivan70,sullivan74}\cite[VI.8.1]{bk} and the theory of $p$-compact
groups \cite{dwyer98,DW94,AGMV03} to translate the
conjecture into a purely combinatorial statement. We then proceed to
show that this statement is `generically' false, with a
counterexample appearing in rank $66$. On the other hand we verify
using a computer  that our counterexample is in fact of minimal rank,
i.e., that the conjecture is true for any finite loop space of rank
less than $66$. This extends earlier work of many authors which show
the statement to be true when the rank is at most $5$ (see
\cite{ochiai68, smith68, sugawara71,hubbuck71,ewing72} and also
\cite{aguade81, lin88, LW91, lin92}).

We now explain this in more detail. Recall that a finite connected
loop space is a triple $(Y,BY,e)$ where $Y$ is a finite connected
CW-complex, $BY$ is a based CW-complex, and $e: Y \to \Omega BY$ is a
homotopy equivalence, where $\Omega BY$ denotes the space of based
loops in $BY$. (We usually refer to a loop space just as $Y$
suppressing the rest of the structure.) It is an old theorem of Hopf
\cite[Satz I]{hopf41} that the rational cohomology of any connected
finite loop space is a graded exterior algebra $H^*(Y;\Q) \cong
\bigwedge_\Q(x_1,...,x_r)$, where the generator $x_i$ is in odd
dimension $2d_i-1$. The number $r$ is called the {\em rank} of $Y$ and the
collection of $d_i$'s are called the {\em degrees} (or if doubled the
`type') of $Y$. It is a classical result of Serre \cite{serre53} that the
collection $\{d_1,\ldots,d_r\}$ in fact uniquely determines the
rational homotopy type of $(Y,BY,e)$. The $p$-completion
$(Y\pcom,BY\pcom,e\pcom)$ is a connected $p$-compact group, i.e.,
$H^*(Y\pcom;\F_p)$ is finite dimensional and connected and $BY\pcom$
is $p$-complete (in the sense of Sullivan \cite{sullivan70,sullivan74}
or Bousfield-Kan \cite{bk}; see e.g., \cite{DW94}, \cite{dwyer98}, or
\cite{AGMV03} for much more on $p$-compact groups.)

An amazing result of Dwyer-Wilkerson \cite{DW94}, extending work of
Dwyer-Miller-Wilkerson \cite{DMW92} and Adams-Wilkerson \cite{AW80},
says that any $p$-compact group $X$ has a maximal torus, that is a
loop map $T \cong (S^1\pcom)^r \to X$ which is suitably maximal and an
associated Weyl group $W_X$.  If $X$ is connected, then $W_X$ acts
faithfully on $L_X = \pi_1(T)$, in such a way that $(W_X,L_X)$ becomes
a finite $\Z_p$-reflection group and
$$H^*(BX;\Z_p) \otimes \Q \pil{\cong} (H^*(BT;\Z_p)\otimes \Q)^{W_X}.$$
The invariant ring $(H^*(BT;\Z_p) \otimes \Q)^{W_X}$ is a polynomial
algebra with generators in dimensions $2e_1, \ldots, 2e_r$, and the
integers $e_1, \ldots, e_r$ are just the well known degrees of the
$\Q_p$-reflection group $(W_X,L_X \otimes \Q)$
\cite[Ch.~7]{benson93},\cite{GM03}. (In fact the harder classification
states that $(W_X,L_X)$ completely classifies $X$ when $p$ is odd
\cite{AGMV03}.) If $Y$ is a finite loop space then, for all primes $p$,
$$H^*(BY;\Q) \otimes \Q_p \cong H^*(BY;\Z) \otimes \Q_p \cong
H^*(BY\pcom;\Z_p) \otimes \Q.$$
Furthermore by the Eilenberg-Moore spectral sequence $H^*(BY;\Q)$ is a
polynomial algebra with generators in dimensions $2d_1, \ldots ,2d_r$,
where the $d_i$'s are the degrees of $Y$ introduced earlier. Hence by
the above we conclude that, for each prime $p$, $Y\pcom$ has to be
some $p$-compact group such that the degrees of $Y$ match up with the
degrees of the $\Q_p$-reflection group $(W_{Y\pcom},L_{Y\pcom} \otimes
\Q)$. This puts severe restrictions on the possible degrees. Finite
$\Q_p$-reflection groups have been classified by Clark-Ewing
\cite{CE74} building on the classification over $\C$ by Shephard-Todd
\cite{ST54}. The classification divides into three infinite families
along with $34$ sporadic cases (the non-Lie ones only being realizable
for certain primes). We denote by $W_{D}$ the Weyl group coming from
the Dynkin diagram $D$ whereas the notation $G(\cdot,\cdot,\cdot)$
means a given group from the infinite family $2$, and $G_n$ refers to
one of the other exotic cases, in the standard notation listed e.g.,
in \cite{CE74} or \cite{GM03}.  Historically, pioneering work of
Clark  \cite{clark63} had already shown that if the maximal degree of
$Y$ is $h$ then $Y$ also has to have the degree $m \leq h$ if $m-1$
and $h$ are relatively prime, using arguments only involving large
primes (compare also \cite[3.20]{humphreys90}). This, as input to
small rank calculations, served as original motivation for the
conjecture. Adams-Wilkerson \cite{AW80} much later found the
restrictions imposed by reflection groups described above,
but worked only at large primes since the technology of \cite{DW94}
was not available. They furthermore gave an example
\cite[Ex.~1.4]{AW80} showing that the large prime information is
algebraically not enough to settle the conjecture. What we
show here is that, contrary to general expectation, the restrictions
at all primes are not even sufficient.

\begin{thm} \label{rank66}
There exists a connected finite loop space $Y$ of rank $66$
such that $Y\pcom$ has Weyl group as a $\Q_p$-reflection group given
by
\begin{alignat*}{3}
W_{A_4B_4B_5B_8B_8E_8A_{12}B_{14}} &\times G_{24} && \qquad \mbox{ for } p=2,\\
W_{A_4D_4B_5B_8B_8E_8A_{13}B_{14}} &\times G_{12} && \qquad \mbox{ for } p \equiv 1,3 \pmod{8},\\
W_{G_2A_4B_4B_4B_7D_{10}A_{13}B_{15}} &\times G(4,2,7) && \qquad \mbox{ for } p \equiv 5 \pmod{8},\\
W_{A_4B_4B_5B_8B_8E_8A_{13}B_{14}} &\times G(6,3,2) && \qquad \mbox{ for } p \equiv 7 \pmod{24},\\
W_{D_4A_5D_8D_8B_{10}A_{13}D_{16}} &\times G(24,24,2) && \qquad \mbox{ for } p \equiv 23 \pmod{24}.
\end{alignat*}
The space $Y$ has dimension $1254$ and the degrees of $Y$ are
$$\{2^8,3^2,4^8,5^2,6^7,7,8^7,9,10^5,11,12^5,13,14^5,16^3,18^2,20^2,22,24^2,26,28,30\}$$
using exponent notation to denote repeated degrees, and these do not
agree with the degrees of any $\Q$-reflection group, i.e., the graded
vector space $H^*(X;\Q)$ does not agree with $H^*(G;\Q)$ for any
compact Lie group $G$.

Furthermore this counterexample is minimal in the sense that any connected finite loop space of rank less than $66$ is rationally equivalent to some compact Lie group $G$.
\end{thm}

A $p$-compact group realizing each of the above simple non-Lie
rational Weyl groups was constructed by Clark-Ewing \cite{CE74} in the
cases where $p$ does not divide the Weyl group order. The remaining
important small prime cases were constructed by  Zabrodsky
$(G_{12},p=3)$ \cite{zabrodsky84}, Dwyer-Wilkerson $(G_{24}, p=2)$
\cite{DW93}, and Notbohm-Oliver $(G(\cdot,\cdot,\cdot), p \mbox{
  small})$ \cite{notbohm98}. Since $G_{24}$ and $G_{12}$ are the only
finite simple $\Q_p$-reflection groups which do not come from compact Lie
groups, for $p=2$ and $3$ respectively, any counterexample will have
to involve these groups.

We note that by work of Bauer-Kitchloo-Notbohm-Pedersen \cite{BKNP03}
the loop space $Y$ of Theorem~\ref{rank66} is in fact homotopy
equivalent to a compact, smooth, parallelizable manifold.

Sullivan's arithmetic square \cite{sullivan70,sullivan74} reduces the study of
finite loop spaces to the study $p$-compact groups for all primes $p$
with the same degrees together with the well understood concept of
arithmetic square `mixing'. We restrict ourselves here to giving the
following lemma which guarantees that an algebraic counterexample
produces a topological counterexample.

\begin{lemma} \label{numerologyred}
Let $\{d_1, \ldots, d_r\}$ be a collection of positive integers (with
repetitions allowed). Suppose that for each prime $p$ we have a
connected $p$-compact group $X_p$ whose Weyl group $(W_{X_p},L_{X_p})$
has degrees $\{d_1, \ldots, d_r\}$. Then there exists a (non-unique)
connected finite loop space $Y$ such that $Y\pcom \cong X_p$ as
$p$-compact groups.
\end{lemma}

The next theorem guarantees that `generically' there will be sets of
degrees which are the degrees of a $\Q_p$-reflection group for all
primes $p$ without being the degrees of any $\Q$-reflection
group. Together with Lemma~\ref{numerologyred} this shows why examples
like the one in Theorem~\ref{rank66} exist. For the statement of the
result we need to introduce some more notation. If $\{d_1, d_2,
\ldots, d_r\}$ is a collection of degrees, then the associated
{\em degree vector} equals $(x_1, x_2, \ldots )\in\Z^{(\infty)}$
where $x_i$ is the number of degrees equal to $i$. We let
$K_{\Lie}\subseteq \Q^{(\infty)}$ denote the positive rational cone
spanned by the degree vectors of the simple $\Q$-reflection groups,
i.e. the set of finite nonnegative rational linear combinations of
these vectors. Similarly we let $K_p$ denote the positive rational
cone spanned by the degree vectors of the simple $\Q_p$-reflection
groups. Finally we let $K_{\Lin}$ denote the positive rational cone
spanned by the degree vectors of the simple $\Q$-reflection groups and
the degree vectors of the groups $G(m,m,2)$, $m=8$, $12$, $24$,
cf. \cite[Thm.~1.1(b)]{lin88}.

\begin{thm}\label{conethm}
We have $K_{\Lie}\subsetneqq \bigcap_p K_p$, where the intersection is
taken over all primes $p$. Moreover $\bigcap_p K_p = K_2
\cap K_3 \cap K_5 \cap K_7 \cap K_{\Lin}$.
\end{thm}

\subsection*{Acknowledgments} We are grateful to Ib Madsen whose
encouraging questions about integrality results in light of the
$p$-complete theory led us to consider this conjecture. We thank
J{\o}rgen Tind for pointing us to Komei Fukuda's program cdd+
\cite{fukuda03} for doing computations on polyhedral cones. The
hospitality of University of Copenhagen, SFB 478 M\"unster, and
\AA rhus University helped make this collaboration possible.

\section{Proofs}
Recall the following essentially classical result.
\begin{thm} \label{sphereresult}
Let $X$ be a simply connected $p$-compact group with degrees $\{d_1,
\ldots, d_r\}$. Then, as a space, $X \simeq  (S^{2d_1-1} \times \cdots
\times S^{2d_r-1})\pcom$ if and only if $p \geq \max\{d_1, \ldots,
d_r\}$. If $X$ is just assumed connected then under the stronger
assumption $p>\max\{d_1, \ldots, d_r\}$, $X$ still splits as a product
of spheres.
\end{thm}

\begin{proof}[Sketch of proof]
Set $h=\max\{d_1, \ldots, d_r\}$ and suppose first that $X$ is simply
connected. If $p \geq h$ then by a Bockstein spectral sequence
argument of Browder \cite[Thm.~4.7]{browder63}, $H^*(X;\Z_p)$ is
torsion free and concentrated in odd degrees. Hence an easy argument of
Serre \cite[Ch.~V~Prop.~6]{serre53} (see also \cite{kumpel72}), using that $\pi_n(S^{2d_i-1})$ has no $p$-torsion when $n<2d_i-1+2p-3$,
yields that $X  \simeq (S^{2d_1-1} \times \cdots \times
S^{2d_r-1})\pcom$. (In fact, this direction uses only that $X$ is an
$H$-space with $H^*(X;\F_p)$ finite dimensional.) The other direction,
which is more subtle and not needed here, was first established for
compact Lie groups by Serre \cite{serre53} and Kumpel \cite{kumpel65}
by case-by-case arguments, and later a general argument was given by
Wilkerson \cite[Thm.~4.1]{wilkerson73} using operations in
$K$-theory.

Assume now that $X$ is just connected and that $p > h$. If $\pi_1(X)$ is
torsion free, then as a space $X \simeq \widetilde X \times
(S^1\pcom)^k$ where $\widetilde X$ is simply connected (cf.~e.g.~\cite[p.~24]{kane88}), which reduces us to the previous case. Hence we just
have to justify that with $p$ as above $\pi_1(X)$ does not have
torsion. By \cite[Thm.~1.4]{MN94} we have a fibration $BK \to
B\widetilde X \times BT' \to BX$ such that $\widetilde X$ is a simply
connected $p$-compact group, $T'$ is a torus, $K$ is a finite $p$-group,
and the projection map $K \to \widetilde X$ is a central
monomorphism. But then by \cite[Thm.~7.6]{dw:center} $K$ is contained
in $\dT^{W_{\widetilde X}}$, where $\dT$ is a discrete approximation
to maximal torus in $\widetilde X$. In particular if
$\dT^{W_{\widetilde X}} = 0$, $\pi_1(X)$ has to be torsion free. But
if $p>h$ then in particular $p \nmid |W_{\widetilde X}|$, so we have
an exact sequence
$$\cdots \to H^0(W_{\widetilde X};L_{\widetilde X}\otimes \Q) \to
H^0(W_{\widetilde X};\dT) \to H^1(W_{\widetilde X};L_{\widetilde X})
\to \cdots$$
where the first and third terms are zero so  $\dT^{W_{\widetilde X}} =
0$ as wanted.
\end{proof}

\begin{proof}[Proof of Lemma~\ref{numerologyred}]
Set $h = \max\{d_1,\ldots, d_r\}$ and let $BM = (\prod_p BX_p)_\Q$. Since
rationalization commutes with taking loop space and finite products we
have that $\Omega BM \simeq
(\prod_{p<h} (X_p)_\Q) \times (\prod_{p\geq h} X_p)_\Q$.
By Theorem~\ref{sphereresult}  $X_p \simeq (S^{2d_1-1} \times \cdots
\times S^{2d_r-1})\pcom$ when $p>h$ and by the same argument $(X_p)_\Q
\simeq ((S^{2d_1-1} \times \cdots \times S^{2d_r-1})\pcom)_\Q$ for all
primes $p$.
Combined with the fact that $(S^{2d_i-1}\pcom)_\Q \simeq
K(\Q_p,2d_i-1)$ this implies that $\pi_n(BM) = (\prod_n
\pi_n(BX_p))_\Q \cong \bigoplus_{i, 2d_i =n} \A_f$, where
$\A_f = (\prod_p \Z_p) \otimes \Q$ is the ring of finite adeles. In
particular $BM$ only has homotopy groups in even dimensions.  But now
a rational space which only has homotopy groups in even dimensions is
necessarily a product of Eilenberg-Mac Lane spaces, as is easily seen
by going up the Postnikov tower (cf. e.g., \cite[Ch.~IX]{whitehead78}). Set $BK = K(\Q,2d_1) \times \cdots
K(\Q,2d_r)$ and construct a map $BK \to BM$ by levelwise taking the
unit ring map $\Q \to \A_f$.
Define $BY$ as the homotopy pullback of the diagram $BK \to BM
\leftarrow \prod_p BX_p$.

Since $\Q$ and $\widehat \Z = \prod_p \Z_p$ generate $\A_f$ the
Mayer-Vietoris sequence in homotopy groups corresponding to a
homotopy pull-back in fact splits, so $\pi_n(BY)$ is the pull-back
in groups of of the diagram $\pi_n(BK) \to \pi_n(BM) \leftarrow
\prod_p \pi_n(BX_p)$. Concretely, the homotopy groups of $BY$ are
given by
$$\pi_n(BY) = (\bigoplus_{i,2d_i=n}\Z) \oplus (\bigoplus_p
\Tor(\Z,\pi_n(BX_p))).$$
In particular, this shows that $\pi_{n}(\Omega BY)$ is finitely
generated for all $n$. Hence also $H^n(\Omega BY;\Z)$ is finitely
generated for all $n$ (see \cite{serre53}\cite[Thm.~2.16]{HMR75}).

By construction $H^*(BY;\Z_p) \pil{\cong} H^*(BX_p;\Z_p)$ so, since the
spaces involved are simply connected, $H^*(\Omega BY;\Z_p) \pil{\cong}
H^*(X_p;\Z_p)$ for all $p$. But since we have seen that each
$H^n(\Omega BY;\Z_p)$ is finitely generated and we know that $X_p$ is
homotopy equivalent to a product of spheres for all $p>h$, we conclude
that in fact $\bigoplus_n H^n(\Omega BY;\Z)$ is finitely generated as an
abelian group.

If $\Omega BY$ is simply connected then it follows from the classical
results of Wall \cite[Thm.~B+F]{wall65} that
$\Omega BY$ is weakly homotopy equivalent to a finite CW-complex
$Y$. If $\Omega BY$ is not simply connected then the conclusion still
holds, now appealing to a more recent result of Notbohm
\cite{notbohm99wall} (see also \cite{BKNP03}) which relies on $\Omega
BY$ being a loop space.
\end{proof}

\begin{proof}[Proof of Theorem~\ref{rank66}]
It follows directly from Lemma~\ref{numerologyred} that we can
construct a connected finite loop space $Y$ with the listed
properties. One can check directly by a finite search that the degrees
of $Y$ does not agree with those of a compact Lie group, but one can also
argue more simply as follows. There are exactly $61$ simple
$\Q$-reflection groups whose degrees are all at most $30$. The inner
product of the vector
\begin{equation} \label{notlietype}
(0,2,-1,-1,0,0,1,-1,0,0,0,1,-1,0,1,-1,0,0,0,0,0,4,-1,-3,1,-1,0,0,0,1)
\end{equation}
and the degree vector of any of these is non-negative. Hence the same
holds for any $\Q$-reflection groups whose degrees are all at most
$30$. However the inner product of the vector \eqref{notlietype} with
the degree vector of $Y$ equals $-1$, so the degrees of $Y$ does not
agree with those of a compact Lie group.

To check that our counterexample has minimal rank we proceed as
follows. For any prime $p$, there are only finitely many
$\Q_p$-reflection groups of a given rank, cf.~\cite{CE74}. Hence one
can go through the list of (say) $\Q_2$-reflection groups and check
which of these has degrees not matching those of any $\Q$-reflection
group, but matching those of a $\Q_p$-reflection group for $p=3, 5,
\ldots$. We have written a C++ program which implements this
algorithm, and used it to check all $\Q_2$-reflection groups of rank
less than $66$.
\end{proof}

\begin{rem}
Note that Theorem~\ref{rank66} in particular tells us that there
exists a loop space whose $p$-completion is not homotopy equivalent to
the $p$-completion of a compact Lie group for any prime $p$. If we
only want this to hold for a single prime $p$ we can find much simpler
examples.  For instance using Lemma~\ref{numerologyred} one can
construct a finite loop space which $2$-completed is homotopy
equivalent to $X_2 = DI(4) \times \Sp(1)\twocom \times \Sp(6)\twocom$
and which $p$-completed for $p\neq 2$ is homotopy equivalent to $X_p =
\Sp(3)\pcom \times \Sp(7)\pcom$. However, $X_2$ is not homotopy
equivalent to the $2$-completion of a compact Lie group, since the
only Lie groups with the right rational degrees are quotients of
$\Sp(3) \times \Sp(7)$, but these do not have the same mod $2$
Poincar\'e series as $X_2$, as can be obtained from \cite{DW93}. (Compare \cite[Conj.~2]{wilkerson74},
\cite[Prob.~9]{stasheff71}, \cite[Conj.~B+C]{kane88}.)
\end{rem}

\begin{proof}[Proof of Theorem~\ref{conethm}]
The first claim follows from the proof of Theorem~\ref{rank66}, see
also Remark~\ref{conestrategy} below. To show the second claim note that the proof of \cite[Thm.~1.1(b)]{lin88} shows that $\bigcap_{p} K_p
\subseteq K_{\Lin}$ and hence
$$\bigcap_{p} K_p \subseteq K_2\cap K_3\cap K_5\cap K_7\cap
K_{\Lin}.$$
To prove the reverse inclusion, note that the only simple $\Q_3$-reflection
group which is not a $\Q$-reflection group is the group
$G_{12}$. Since this is a $\Q_p$-reflection group for all primes $p$
satisfying $p\equiv 1,3 \pmod{8}$ we get $K_3\subseteq K_p$ for these
primes. Similarly the simple $\Q_5$-reflection groups which are not
$\Q$-reflection groups are precisely the groups
$G(4,1,n), G(4,2,n)$ for $n\geq 2$, $G(4,4,n)$ for $n\geq 3$, $\Z/4$
from family $3$, $G_8$, $G_{29}$ and $G_{31}$. This shows that
$K_5\subseteq K_p$ when $p\equiv 1 \pmod{4}$. In the same way we see
that $K_7\subseteq K_p$ when $p$ satisfies $p\equiv 1 \pmod{6}$ and
$p\equiv \pm 1 \pmod{8}$, i.e. when $p\equiv 1,7 \pmod{24}$. Finally
the groups $G(m,m,2)$, $m=8$, $12$, $24$ are all $\Q_p$-reflection
groups when $p\equiv \pm 1 \pmod{24}$, so $K_{\Lin}\subseteq K_p$ for
these primes. This proves the result
since any prime $p$ satisfies $p=2$, $p\equiv 1,3 \pmod{8}$, $p\equiv
1 \pmod{4}$, $p\equiv 1,7 \pmod{24}$ or $p\equiv \pm 1 \pmod{24}$.
\end{proof}

\begin{rem} \label{conestrategy}
A few remarks about how we found the counterexample in
Theorem~\ref{rank66} might be in order, since this is not really clear
from the proof. First we used Fukuda's cdd+ program \cite{fukuda03} to
establish Theorem~\ref{conethm} by showing that $K_{\Lie}$ truncated
at degree say $30$ does not agree with the intersections of the
similarly truncated versions of $K_2$, $K_3$, $K_5$, $K_7$, and
$K_{\Lin}$. From this it is a linear programming problem
to obtain a concrete point in the difference, and by solving the
associated integer programming problem one gets a point in the
difference which is minimal with respect to rank (or dimension). Note
however that being minimal in this sense is slightly weaker
than being a minimal counterexample, which is why we had to finish off
our proof of minimality in Theorem~\ref{rank66} with a brute force
check, which required rather massive computer calculations. Using the
geometric picture we have found a counterexample of smaller dimension
but larger rank. Namely, there exists a connected finite loop space
with rank $74$ and dimension $1250$ and degrees
$$\{2^9,3^2,4^7,5^3,6^8,7^3,8^8,9^3,10^6,11^2,12^6,13^2,14^5,15,16^3,18^2,20^2,22,24\},$$
which is not rationally equivalent to a compact Lie group. This is
seen by considering the $\Q_p$-reflection groups
\begin{alignat*}{3}
W_{A_1A_1E_6D_7D_8A_9D_{11}D_{13}A_{15}} &\times G_{24} && \qquad \mbox{ for } p=2,\\
W_{A_1A_1D_5D_7D_9B_{10}B_{12}A_{13}A_{14}} &\times G_{12} && \qquad \mbox{ for } p\equiv 1,3 \pmod{8},\\
W_{A_1A_1D_5D_6E_7D_9B_{11}A_{13}A_{14}} &\times G(4,4,7) && \qquad
\mbox{ for } p\equiv 5 \pmod{8},\\
W_{G_2D_5D_7D_9B_{10}B_{12}A_{13}A_{14}} &\times G(8,8,2) && \qquad
\mbox{ for } p\equiv 7 \pmod{8}.
\end{alignat*}
It is also possible to construct a counterexample where all degrees
are even, cf. \cite{LW91}. For instance there is an example of rank
$68$, dimension $1468$ and degrees
$$\{2^8,4^8,6^8,8^8,10^6,12^5,14^6,16^4,18^3,20^3,22^2,24^3,26^2,28,30\}.$$
Here one can use the $\Q_p$-reflection groups
\begin{alignat*}{3}
W_{B_4B_5B_5B_8B_8E_8B_{13}B_{14}} &\times G_{24} && \qquad \mbox{ for } p=2,\\
W_{D_4B_5B_5B_8B_8E_8B_{14}D_{14}} &\times G_{12} && \qquad \mbox{ for } p \equiv 1,3 \pmod{8},\\
W_{G_2B_4B_5B_5B_7B_9D_{14}B_{15}} &\times G(4,2,7) && \qquad \mbox{ for } p
\equiv 5 \pmod{8},\\
W_{B_4B_5B_5B_8B_8E_8B_{14}D_{14}} &\times G(6,3,2) && \qquad \mbox{ for } p \equiv 7 \pmod{24},\\ 
W_{D_4D_6D_8D_8B_{10}D_{14}D_{16}} &\times G(24,24,2) && \qquad \mbox{ for } p
\equiv 23 \pmod{24}.
\end{alignat*}
\end{rem}

\bibliographystyle{plain}
\bibliography{poddclassification}
\end{document}